\documentclass[12pt,leqno]{article}
\usepackage[doc]{optional}
\usepackage{color}
\usepackage{soul}
\usepackage{graphicx}
\definecolor{labelkey}{rgb}{0,0.08,0.45}
\definecolor{refkey}{rgb}{0,0.6,0.0}
\definecolor{Brown}{rgb}{0.45,0.0,0.05}
\definecolor{lime}{rgb}{0.00,0.8,0.0}
\definecolor{lblue}{rgb}{0.5,0.5,0.99}
\usepackage{mathpazo}
\usepackage{amsmath}
\usepackage{amssymb}
\usepackage{theorem}
\newcommand{\weakly}{\ensuremath{\:{\rightharpoonup}\:}}

\newcommand{\nnn}{\ensuremath{{n\in{\mathbb N}}}}

\newcommand{\menge}[2]{\big\{{#1}~\big |~{#2}\big\}}
\newcommand{\mmenge}[2]{\bigg\{{#1}~\bigg |~{#2}\bigg\}}

\newcommand{\To}{\ensuremath{\rightrightarrows}}

\newcommand{\fenv}[1]%
{\ensuremath{\,\overrightarrow{\operatorname{env}}_{#1}}}
\newcommand{\benv}[1]%
{\ensuremath{\,\overleftarrow{\operatorname{env}}_{#1}}}

\newcommand{\scal}[2]{\left\langle{#1},{#2}  \right\rangle}

\newcommand{\zeroun}{\ensuremath{\left]0,1\right[}}
\newcommand{\RR}{\ensuremath{\mathbb R}}
\newcommand{\RP}{\ensuremath{\left[0,+\infty\right[}}
\newcommand{\RPX}{\ensuremath{\left[0,+\infty\right]}}
\newcommand{\RPP}{\ensuremath{\,\left]0,+\infty\right[}}
\newcommand{\RX}{\ensuremath{\,\left]-\infty,+\infty\right]}}

\newcommand{\dom}{\ensuremath{\operatorname{dom}}}

\newcommand{\gr}{\ensuremath{\operatorname{gr}}}
\newcommand{\prox}{\ensuremath{\operatorname{Prox}}}

\newcommand{\ran}{\ensuremath{\operatorname{ran}}}

\newcommand{\Fix}{\ensuremath{\operatorname{Fix}}}
\newcommand{\Id}{\ensuremath{\operatorname{Id}}}

\newcommand{\minf}{\ensuremath{-\infty}}
\newcommand{\pinf}{\ensuremath{+\infty}}

{\begin{list}{}{%
\settowidth{\labelwidth}{\textrm{#1~}}%
\setlength{\leftmargin}{\labelwidth+\labelsep}}}
{\end{list}}
\newtheorem{theorem}{Theorem}[section]
\newtheorem{lemma}[theorem]{Lemma}
\newtheorem{corollary}[theorem]{Corollary}
\newtheorem{proposition}[theorem]{Proposition}
\newtheorem{definition}[theorem]{Definition}
\theoremstyle{plain}{\theorembodyfont{\rmfamily}
}
\theoremstyle{plain}{\theorembodyfont{\rmfamily}
}
\theoremstyle{plain}{\theorembodyfont{\rmfamily}
}
\theoremstyle{plain}{\theorembodyfont{\rmfamily}
\newtheorem{example}[theorem]{Example}}
\newtheorem{fact}[theorem]{Fact}
\theoremstyle{plain}{\theorembodyfont{\rmfamily}
\newtheorem{remark}[theorem]{Remark}}

\newcommand{\boxedeqn}[1]{%
    \[\fbox{%
        \addtolength{\linewidth}{-2\fboxsep}%
        \addtolength{\linewidth}{-2\fboxrule}%
        \begin{minipage}{\linewidth}%
        \begin{equation}#1\\[+4mm]\end{equation}%
        \end{minipage}%
      }\]%
  }

\begin{document}

\title{\textrm{
Firmly nonexpansive mappings\\ and\\ maximally monotone operators:\\
correspondence and duality}}

\author{
Heinz H.\ Bauschke\thanks{Mathematics, Irving K.\ Barber School,
UBC Okanagan, Kelowna, British Columbia V1V 1V7, Canada. E-mail:
\texttt{heinz.bauschke@ubc.ca}.},
~Sarah M.\ Moffat\thanks{Mathematics, Irving K.\ Barber School, UBC Okanagan,
Kelowna, British Columbia V1V 1V7, Canada. E-mail:
\texttt{sarah.moffat@ubc.ca}.},
~and Xianfu\
Wang\thanks{Mathematics, Irving K.\ Barber School, UBC Okanagan,
Kelowna, British Columbia V1V 1V7, Canada.
E-mail:  \texttt{shawn.wang@ubc.ca}.}}

\date{January 24, 2011}
\maketitle

\vskip 8mm

\begin{abstract} \noindent
The notion of a firmly nonexpansive mapping is central in fixed point
theory because of attractive convergence properties for iterates and the
correspondence with maximal monotone operators due to Minty. 
In this paper, we systematically analyze the relationship between
properties of firmly nonexpansive mappings and associated maximal monotone
operators.
Dual and self-dual properties are also identified. 
The results are illustrated through several examples. 

\end{abstract}

{\small


}

\section{Introduction}

Throughout this paper,
\boxedeqn{
\text{$X$ is
a real Hilbert space with inner
product $\scal{\cdot}{\cdot}$ and induced norm $\|\cdot\|$.}
}
Recall that a mapping
\begin{equation}
T\colon X\to X
\end{equation}
is \textbf{firmly nonexpansive}
if 
\begin{equation}
(\forall x\in X)(\forall y\in X)
\quad
\|Tx-Ty\|^2 + \|(\Id-T)x-(\Id-T)y\|^2 \leq \|x-y\|^2,
\end{equation}
where $\Id\colon X\to X\colon x\mapsto x$ denotes
the \textbf{identity operator}. 
It is clear that if $T$ is firmly nonexpansive,
then it is \textbf{nonexpansive}, i.e., 
{Lipschitz continuous} with constant $1$, 
\begin{equation}
(\forall x\in X)(\forall y\in X) \quad
\|Tx-Ty\|\leq\|x-y\|;
\end{equation}
the converse, however, is false (consider $-\Id$).
When $T$ is Lipschitz continuous with a constant in 
$\left[0,1\right[$, then we shall refer to $T$ as a
\textbf{Banach contraction}. 

Returning to firmly nonexpansive mappings, 
let us provide a useful list of well-known characterizations.

\begin{fact}
\label{f:firm}
{\rm (See, e.g., \cite{BC2011,GK,GR}.)}
Let $T\colon X\to X$.
Then the following are equivalent:
\begin{enumerate}
\item $T$ is firmly nonexpansive.
\item $\Id-T$ is firmly nonexpansive.
\item $2T-\Id$ is nonexpansive.
\item $(\forall x\in X)(\forall y\in X)$
$\|Tx-Ty\|^2 \leq \scal{x-y}{Tx-Ty}$. 
\item 
\label{f:firmsymm}
$(\forall x\in X)(\forall y\in X)$.
$0\leq \scal{Tx-Ty}{(\Id-T)x-(\Id-T)y}$.
\end{enumerate}
\end{fact}

Various problems in the natural sciences and engineering
can be converted into a fixed point problem, where the set of
desired solutions is the set of \textbf{fixed points}
\begin{equation}
\Fix T := \menge{x\in X}{x=Tx}.
\end{equation}
If $T$ is firmly nonexpansive and $\Fix T\neq\varnothing$,
then the sequence of iterates
\begin{equation}
(T^nx)_\nnn
\end{equation}
is well known to converge weakly to a fixed point \cite{Browder67} ---
this is not true for mappings that are merely nonexpansive; consider,
e.g., $-\Id$.

Firmly nonexpansive mappings are also important because of their
correspondence with maximally monotone operators.
Recall that a set-valued operator $A\colon X\To X$ (i.e., $(\forall x\in
X)$ $Ax\subseteq X$) with {graph} $\gr A$ is \textbf{monotone} if
\begin{equation}
(\forall (x,u)\in\gr A)
(\forall (y,v)\in\gr A)
\quad
\scal{x-y}{u-v}\geq 0,
\end{equation}
and that $A$ is \textbf{maximally monotone} if it is impossible to find a
proper extension of $A$ that is still monotone.
We write $\dom A :=\menge{x\in X}{Ax\neq\varnothing}$
and $\ran A := A(X) = \bigcup_{x\in X}Ax$ for the
\textbf{domain} and \textbf{range} of $A$, respectively. 
Monotone operators play a crucial role in modern nonlinear analysis and
optimization; see, e.g., the books
\cite{BC2011}, 
\cite{BorVanBook}, 
\cite{Brezis}, 
\cite{BurIus},
\cite{Simons1},
\cite{Simons2},
\cite{Zalinescu}, 
\cite{Zeidler2a},
\cite{Zeidler2b}, 
and 
\cite{Zeidler1}.

Prime examples of maximally monotone operators are
continuous linear monotone operators and subdifferential operators 
(in the sense of convex analysis) of 
functions that are convex, lower semicontinuous, and proper. It will be
convenient to set 
\boxedeqn{
\Gamma_0 := \menge{f\colon X\to\RX}{\text{$f$ is convex, lower
semicontinuous, and proper}}. 
}
See  \cite{Rock70}, \cite{Rock98}, \cite{Zalinescu}, and 
\cite{BC2011} for background material on convex analysis. 

Now let $A\colon X\To X$ be maximally monotone
and denote the associated \textbf{resolvent} by 
\begin{equation}
J_A := (\Id+A)^{-1}.
\end{equation}
When $A=\partial f$ for some $f\in\Gamma_0$, i.e., $A$ is a 
\textbf{subdifferential operator}, then we also write $J_{\partial f} =
\prox_f$ and, following Moreau \cite{Moreau},
we refer to this mapping as the \textbf{proximal mapping}. 
In his seminal paper \cite{Minty}, Minty observed 
that $J_A$ is in fact a firmly nonexpansive operator
from $X$ to $X$ and that, conversely, every firmly nonexpansive operator
arises this way:
\begin{fact}[Minty] 
{\rm (See, e.g., \cite{Minty} or \cite{EckBer}.)}
Let $T\colon X\to X$ be firmly nonexpansive, 
and let $A\colon X\To X$ be maximally monotone. 
Then the following hold.
\begin{enumerate}
\item 
$B= T^{-1}-\Id$ is maximally monotone (and $J_B = T$).
\item
$J_A$ is firmly nonexpansive (and $A = J_A^{-1}-\Id$).
\end{enumerate}
\end{fact}

Therefore, the mapping 
\begin{equation}
T \mapsto T^{-1}-\Id
\end{equation}
from the set of firmly nonexpansive mappings to the
set of maximal monotone operators
is a bijection, with inverse
\begin{equation}
A \mapsto J_A = (\Id + A)^{-1}.
\end{equation}

Each class of objects has its own duality operation:
if $T\colon X\to X$ is firmly nonexpansive,
then so is the associated \textbf{dual mapping} (see Fact~\ref{f:firm})
\begin{equation}
\Id-T;
\end{equation}
note that the dual of the dual is indeed $\Id-(\Id-T) = T$.
Analogously, 
if $A\colon X\To X$ is maximally monotone,
then so is the inverse operator
\begin{equation}
A^{-1},
\end{equation}
and we clearly have $(A^{-1})^{-1} = A$. 
Traversing between the two classes and dualizing is commutative
as the 
\textbf{resolvent identity}
\begin{equation}
\label{e:resident}
\Id = J_A + J_{A^{-1}}
\end{equation}
shows. 
We also have the \textbf{Minty parametrization}
\begin{equation}
\label{e:Mintypar}
\gr A = \menge{\big(J_Ax,(\Id-J_A)x\big)}{x\in X},
\end{equation}
which provides the bijection $x\mapsto (J_Ax,x-J_Ax)$ 
from $X$ onto $\gr A$, with inverse $(x,u)\mapsto x+u$. 

{\em 
The goal of this paper is threefold.
First, we will provide a comprehensive catalogue of corresponding
properties of firmly nonexpansive mappings and 
maximally monotone operators.
Second, we shall examine duality of these properties, 
and also identify those properties
that are self-dual, i.e., a property holds for $T$ if and only if the
same property holds for $\Id-T$, and a corresponding property holds
for $A$ if and only if it holds for $A^{-1}$ in the context of 
maximally monotone operators.
Third, we revisit some of these properties in more specialized
settings and present applications of our results to operators
occurring in splitting methods, including reflected resolvents.
}

The paper is organized as follows.
In Section~\ref{s:2}, we systematically study properties of
firmly nonexpansive mappings and corresponding properties of
maximally monotone operators.
Section~\ref{s:3} revisits these properties from a duality point of
view. 
In the final Section~\ref{s:4}, we
also touch upon corresponding properties of nonexpansive mappings
and provide a few possible applications to mappings arising in
algorithms. 

The notation we employ is standard and as in 
\cite{BC2011}, \cite{Rock70}, \cite{Rock98},
or \cite{Zalinescu}. These books provide also background
material for notions not explicitly reviewed here.

\section{Correspondence of Properties}

\label{s:2}

The main result in this section is the following result, 
which provides a comprehensive list of corresponding properties
of firmly nonexpansive mappings and maximally monotone
operators.

\begin{theorem}
\label{t:1}
Let $T\colon X\to X$ be firmly nonexpansive,
let $A\colon X\To X$ be maximally monotone, and
suppose that $T = J_A$ or equivalently that $A = T^{-1}-\Id$. 
Then the following hold:
\begin{enumerate}
\item 
\label{t:1i}
$\ran T = \dom A$.
\item
\label{t:1ii}
$T$ is surjective if and only if $\dom A = X$.
\item
\label{t:1ii+}
$\Id-T$ is surjective if and only if $A$ is surjective. 
\item
\label{t:1iii}
$T$ is injective if and only if $A$ is at most single-valued.
\item
\label{t:1isom}
$T$ is an \textbf{isometry}, i.e.,
\begin{equation}
(\forall x\in X)(\forall y\in X)
\quad
\|Tx-Ty\| = \|x-y\|
\end{equation}
if and only if there exists $z\in X$ such that 
$A\colon x\mapsto z$, in which case $T\colon x\mapsto x-z$. 
\item
\label{t:1iv}
$T$ satisfies
\begin{equation}
\label{e:imagestrictfirm}
(\forall x\in X)(\forall y\in X)
\quad
Tx\neq Ty
\;\Rightarrow\;
\|Tx-Ty\|^2 < \scal{x-y}{Tx-Ty}
\end{equation}
if and only if 
$A$ is \textbf{strictly monotone}, i.e.,
\begin{equation}
(\forall (x,u)\in \gr A)(\forall (y,v)\in \gr A)
\quad
x\neq y
\;\Rightarrow\;
\scal{x-y}{u-v}>0.
\end{equation}
\item
\label{t:1iv+}
$T$ is strictly monotone if and only if $A$ is at most single-valued.
\item
\label{t:1v}
$T$ is \textbf{strictly firmly nonexpansive}, i.e., 
\begin{equation}
(\forall x\in X)(\forall y\in X)
\quad
x\neq y
\;\Rightarrow\;
\|Tx-Ty\|^2 < \scal{x-y}{Tx-Ty}
\end{equation}
if and only if 
$A$ is at most single-valued and strictly monotone.
\item 
\label{t:1vi}
$T$ is \textbf{strictly nonexpansive}, i.e.,
\begin{equation}
(\forall x\in X)(\forall y\in X)
\quad
x\neq y
\;\Rightarrow\;
\|Tx-Ty\| < \|x-y\|
\end{equation}
if and only if $A$ is \textbf{disjointly
injective}\footnote{In some of the literature, $A$ is simply referred
to as ``injective'' but this is perhaps a little imprecise: 
``injectivity of $A$'' could be interpreted as 
``[$\{x,y\}\subseteq\dom A$ and $x\neq y$]~$\Rightarrow$~$Ax\neq Ay$'' 
which is different from what we require here.}, 
i.e.,
\begin{equation}
(\forall x\in X)(\forall y\in X)
\quad
x\neq y
\;\Rightarrow\;
Ax\cap Ay = \varnothing. 
\end{equation}
\item
\label{t:1vii}
$T$ is injective and strictly nonexpansive, i.e.,
\begin{equation}
(\forall x\in X)(\forall y\in X)
\quad
x\neq y
\;\Rightarrow\;
0<\|Tx-Ty\| < \|x-y\|
\end{equation}
if and only if $A$ is at most single-valued and disjointly injective. 
\item
\label{t:1viii}
Suppose that $\varepsilon\in\RPP$. 
Then $(1+\varepsilon)T$ is firmly nonexpansive
if and only if 
$A$ is \textbf{strongly monotone} with constant $\varepsilon$, i.e., 
$A-\varepsilon\Id$ is monotone, 
in which case $T$ is
a Banach contraction with constant $(1+\varepsilon)^{-1}$. 
\item
\label{t:1coco}
Suppose that $\gamma\in\RPP$. 
Then $(1+\gamma)(\Id-T)$ is firmly nonexpansive
if and only if $A$ is \textbf{$\gamma$-cocoercive}, i.e.,
\begin{equation}
(\forall (x,u)\in\gr A)(\forall (y,v)\in\gr A)
\quad
\scal{x-y}{u-v}\geq \gamma\|u-v\|^2.
\end{equation}
\item
\label{t:1Banach}
Suppose that $\beta\in\zeroun$. 
Then $T$ is a {Banach contraction} with constant $\beta$
if and only if
$A$ satisfies
\begin{equation}
\label{e:t1Banach}
(\forall (x,u)\in\gr A)(\forall (y,v)\in\gr A)
\quad
\frac{1-\beta^2}{\beta^2}\|x-y\|^2 
\leq 2\scal{x-y}{u-v} + \|u-v\|^2.
\end{equation}
\item
\label{t:1uni}
Suppose that $\phi\colon\RP\to\RPX$ is increasing and vanishes only at
$0$. Then $T$ satisfies
\begin{equation}
(\forall x\in X)(\forall y\in X)
\quad
\scal{Tx-Ty}{(x-Tx)-(y-Ty)}\geq\phi\big(\|Tx-Ty\|\big)
\end{equation}
if and only if
$A$ is \textbf{uniformly monotone} with modulus $\phi$, i.e.,
\begin{equation}
(\forall (x,u)\in\gr A)(\forall (y,v)\in\gr A)
\quad
\scal{x-y}{u-v}\geq\phi\big(\|x-y\|\big).
\end{equation}
\item
\label{t:1para}
$T$ satisfies
\begin{equation}
\label{e:fpara}
(\forall x\in X)(\forall y\in X)
\quad
\|Tx-Ty\|^2 = \scal{x-y}{Tx-Ty}
\;\Rightarrow\;
\begin{cases}
Tx = T\big(Tx+y-Ty\big)\\
Ty = T\big(Ty+x-Tx\big)
\end{cases}
\end{equation}
if and only if $A$ is \textbf{paramonotone} {\rm \cite{CIZ}}, i.e.,
\begin{equation}
(\forall (x,u)\in\gr A)(\forall (y,v)\in\gr A)
\quad
\scal{x-y}{u-v}=0
\;\Rightarrow\;
\{(x,v),(y,u)\}\subseteq \gr A.
\end{equation}
\item
\label{t:1cyc}
{\rm (Bartz et al.)}
$T$ is \textbf{cyclically firmly nonexpansive}, i.e., 
\begin{equation}
\label{e:t1cyc}
\sum_{i=1}^{n}\scal{x_i-Tx_i}{Tx_i-Tx_{i+1}}\geq 0,
\end{equation}
for every set of points $\{x_1,\ldots,x_n\}\subseteq X$,
where $n\in\{2,3,\ldots\}$ and $x_{n+1}=x_1$,
if and only if $A$ is a subdifferential operator,
i.e., there exists $f\in\Gamma_0$ such that $A=\partial f$. 
\item
\label{t:13*}
$T$ satisfies
\begin{equation}
\label{e:t13*}
(\forall x\in X)(y\in X)
\quad
\inf_{z\in X} \scal{Tx-Tz}{(y-Ty)-(z-Tz)} > \minf
\end{equation}
if and only if $A$ is \textbf{rectangular} 
{\rm \cite[Definition~31.5]{Simons2}}
(this is also known as \textbf{$3^*$ monotone}), i.e.,
\begin{equation}
(\forall x\in\dom A)(\forall v\in\ran A)
\quad
\inf_{(z,w)\in\gr A} \scal{x-z}{v-w}>\minf.
\end{equation}
\item
\label{t:1linear}
$T$ is linear if and only if $A$ is a \textbf{linear relation}, i.e.,
$\gr A$ is a linear subspace of $X\times X$. 
\item
\label{t:1affine}
$T$ is affine if and only if $A$ is an \textbf{affine relation}, i.e.,
$\gr A$ is an affine subspace of $X\times X$. 
\item
\label{t:1zara}
{\rm (Zarantonello)}
$\ran T=\Fix T =:C$ if and only if $A$ is a normal cone operator, i.e.,
$A=\partial\iota_C$; 
equivalently, $T$ is a projection (nearest point) mapping $P_C$. 
\item
\label{t:1weak}
$T$ is sequentially weakly continuous
if and only if $\gr A$ is sequentially weakly closed. 
\end{enumerate}
\end{theorem}
\begin{proof}
Let $x,y,u,v$ be in $X$.

\ref{t:1i}:
Clear.

\ref{t:1ii}: 
This follows from \ref{t:1i}. 

\ref{t:1ii+}:
Clear from the Minty parametrization \eqref{e:Mintypar}. 

\ref{t:1iii}: 
Assume first that $T$ is injective
and that $\{u,v\}\subseteq Ax$. 
Then $\{x+u,x+v\}\subseteq (\Id+A)x$
and hence $x=T(x+u)=T(x+v)$.
Since $T$ is injective, it follows that
$x+u=x+v$ and hence that $u=v$. 
Thus, $A$ is at most single-valued.
Conversely,
let us assume that $Tu=Tv=x$. 
Then $\{u,v\} \subseteq (\Id+A)x = x+Ax$
and hence $\{u-x,v-x\}\subseteq Ax$.
Since $A$ is at most single-valued,
we have $u-x=v-x$ and so $u=v$.
Thus, $T$ is injective.

\ref{t:1isom}:
Assume first that $T$ is an isometry.
Then $\|Tx-Ty\|^2 = \|x-y\|^2 \geq \|Tx-Ty\|^2 +
\|(\Id-T)x-(\Id-T)y\|^2$.
It follows that there exists $z\in X$ such that
$T\colon w\mapsto w-z$.
Thus, $T^{-1} \colon w\mapsto w+z$.
On the other hand, $T^{-1}=\Id+A\colon w\mapsto w+Aw$.
Hence $A\colon w\mapsto z$, as claimed.
Conversely, let us assume that
there exits $z\in X$ such that $A\colon w\mapsto z$.
Then $\Id+A\colon w\mapsto w+z$ and hence
$T = J_A = (\Id+A)^{-1}\colon w\mapsto w-z$. 
Thus, $T$ is an isometry. 

\ref{t:1iv}: 
Assume first that $T$ satisfies \eqref{e:imagestrictfirm}, 
that $\{(x,u),(y,v)\}\subseteq \gr A$,
and that $x\neq y$. 
Set $p = x+u$ and $q=y+v$.
Then $(x,u) = (Tp,p-Tp)$ and $(y,v)=(Tq,q-Tq)$. 
Since $x\neq y$, it follows that $Tp\neq Tq$ and
therefore that 
$\|Tp-Tq\|^2 < \scal{p-q}{Tp-Tq}$ 
because $T$ satisfies \eqref{e:imagestrictfirm}. 
Hence $0<\scal{(p-Tp)-(q-Tq)}{Tp-Tq} = \scal{u-v}{x-y}$.
Thus, $A$ is strictly monotone. 
Conversely, let us assume that $A$ is strictly monotone
and that $x=Tu\neq Tv=y$. 
Then $\{(x,u-x),(y,v-y)\}\subseteq \gr A$.
Since $x\neq y$ and $A$ is strictly monotone,
we have $\scal{x-y}{(u-x)-(v-y)}>0$, i.e.,
$\|x-y\|^2 < \scal{x-y}{u-v}$, i.e.,
$\|Tu-Tv\|^2 < \scal{Tv-Tu}{u-v}$.
Thus, $T$ satisfies \eqref{e:imagestrictfirm}. 

\ref{t:1iv+}:
In view of \ref{t:1iv} it suffices to show that
$T$ is injective if and only if $T$ is strictly monotone.
Assume first that $T$ is injective and that $x\neq y$.
Then $Tx\neq Ty$ and hence $0<\|Tx-Ty\|^2 \leq \scal{x-y}{Tx-Ty}$.
Thus, $T$ is strictly monotone. 
Conversely, assume that $T$ is strictly monotone and that $x\neq y$.
Then $\scal{x-y}{Tx-Ty}>0$ and hence $Tx\neq Ty$.
Thus, $T$ is injective.

\ref{t:1v}:
Observe that $T$ is strictly firmly nonexpansive
if and only if $T$ is injective and $T$ satisfies
\eqref{e:imagestrictfirm}. Thus, the result follows
from combining \ref{t:1iii} and \ref{t:1iv}. 

\ref{t:1vi}:
Assume first that $T$ is strictly nonexpansive, 
that $x\neq y$, and to the contrary that
$u\in Ax\cap Ay$. 
Then $x+u\in(\Id+A)x$ and $y+u\in(\Id+A)y$;
equivalently, $T(x+u)=x\neq y = T(y+u)$. 
Since $T$ is strictly nonexpansive, we have
$\|x-y\|=\|T(x+u)-T(y+u)\|<\|(x+u)-(y+u)\|=\|x-y\|$,
which is absurd. Thus, $A$ is disjointly injective. 
Conversely, assume that $A$ is disjointly injective, 
that $u\neq v$, and to the contrary that 
$\|Tu-Tv\|=\|u-v\|$. Since $T$ is firmly nonexpansive,
we deduce that $u-Tu=v-Tv$. 
Assume that $x=u-Tu=v-Tv$.
Then $Tu=u-x$ and $Tv=v-x$;
equivalently, $u-x \in (\Id+A)u$ and
$v-x\in (\Id+A)v$. Thus, $-x\in Au\cap Av$, which contradicts
the assumption on disjoint injectivity of $A$. 

\ref{t:1vii}:
Combine \ref{t:1iii} and \ref{t:1vi}. 

\ref{t:1viii}: 
(See also \cite[Proposition~23.11]{BC2011}.)
Assume first that $(1+\varepsilon)T$ is firmly nonexpansive
and that $\{(x,u),(y,v)\}\subseteq\gr A$. 
Then $x=T(x+u)$ and $y=T(y+v)$.
Hence, $\scal{(x+u)-(y+v)}{x-y}\geq (1+\varepsilon)\|x-y\|^2$
$\Leftrightarrow$ $\scal{x-y}{u-v}\geq\varepsilon\|x-y\|^2$.
Thus, $A-\varepsilon\Id$ is monotone. 
Conversely, assume that $A-\varepsilon\Id$ is monotone
and that $\{(x,u),(y,v)\}\subseteq\gr T$. 
Then $\{(u,x-u),(v,y-v)\}\subseteq\gr A$ and hence
$\scal{u-v}{(x-u)-(y-v)}\geq \varepsilon\|u-v\|^2$;
equivalently,
$\scal{x-y}{u-v}\geq(1+\varepsilon)\|u-v\|^2$. 
Thus, $(1+\varepsilon)T$ is firmly nonexpansive. 

\ref{t:1coco}:
Applying \ref{t:1viii} to $\Id-T$ and $A^{-1}$, we see that
$(1+\gamma)(\Id-T)$ is firmly nonexpansive
if and only if $A^{-1}-\gamma\Id$ is monotone,
which is equivalent to $A$ being $\gamma$-cocoercive.

\ref{t:1Banach}:
Assume first that $T$ is a Banach contraction with constant $\beta$
and that $\{(x,u),(y,v)\}\subseteq\gr A$. 
Set $p=x+u$ and $y=y+v$. 
Then $(x,u) = (Tp,p-Tp)$, $(y,v)=(Tq,q-Tq)$,
and $\|Tp-Tq\|\leq\beta\|p-q\|$, i.e.,
\begin{align}
\|x-y\|^2&\leq\beta^2\|(x+u)-(y+v)\|^2 = \beta^2\|(x-y)+(u-v)\|^2\\
&=\beta^2\big(\|x-y\|^2 + 2\scal{x-y}{u-v}+\|u-v\|^2 \big).\notag
\end{align}
Thus, \eqref{e:t1Banach} holds.
The converse is proved similarly. 

\ref{t:1uni}:
The equivalence is immediate from the Minty parametrization
\eqref{e:Mintypar}. 

\ref{t:1para}:
Assume first that $T$ satisfies \eqref{e:fpara} 
and that $\{(x,u),(y,v)\}\subseteq \gr A$ with
$\scal{x-y}{u-v}=0$. 
Set $p=x+u$ and $q=y+v$.
Then $(x,u)=(Tp,p-Tp)$, $(y,v)=(Tq,q-Tq)$, and
$\scal{Tp-Tq}{(p-Tp)-(q-Tq)} = 0$, i.e.,
$\|Tp-Tq\|^2 = \scal{p-q}{Tp-Tq}$. 
By \eqref{e:fpara}, $Tp=T(Tp+q-Tq)$ and $Tq=T(Tq+p-Tp)$, i.e.,
$x=T(x+v)$ and $y=T(y+u)$, i.e., 
$x+v\in x+Ax$ and $y+u\in y+Ay$, i.e.,
$v\in Ax$ and $u\in Ay$. 
Thus, $A$ is paramonotone. 
Conversely, assume that $A$ is paramonotone,
that $\|Tu-Tv\|^2 = \scal{u-v}{Tu-Tv}$,
that $x=Tu$, and that $y=Tv$.
Then $\{(x,u-x),(y,v-y)\}\subseteq\gr A$
and $\scal{x-y}{(u-x)-(v-y)}=0$. 
Since $A$ is paramonotone, we deduce that
$v-y\in Ax$ and $u-x\in Ay$, i.e.,
$x-y+v\in (\Id+A)x$ and $y-x+u\in (\Id+A)y$, i.e.,
$x=T(x-y+v)$ and $y=T(y-x+u)$, i.e.,
$Tu=T(Tu+v-Tv)$ and $Tv = T(Tv+u-Tu)$. 
Thus, $T$ satisfies \eqref{e:fpara}. 

\ref{t:1cyc}:
This follows from \cite[Theorem~6.6]{BBBRW}. 

\ref{t:13*}:
The equivalence is immediate from the Minty parametrization
\eqref{e:Mintypar}. 

\ref{t:1linear}:
Indeed, 
$T=J_A$ is linear
$\Leftrightarrow$ 
$(A+\Id)^{-1}$ is a linear relation
$\Leftrightarrow$ 
$A+\Id$ is a linear relation
$\Leftrightarrow$ 
$A$ is a linear relation
(see \cite{Cross} for more on linear relations). 

\ref{t:1affine}:
This follows from \ref{t:1linear}. 

\ref{t:1zara}:
This follows from \cite[Corollary~2 on page~251]{Zara}. 

\ref{t:1weak}: 
Assume first that $T$ is sequentially weakly continuous.
Let $(x_n,u_n)_\nnn$ be a sequence in $\gr A$ that converges
weakly to $(x,u)\in X\times X$. 
Then $(x_n+u_n)_\nnn$ converges weakly to 
$x+u$. On the other hand, $\Id-T$ is sequentially weakly continuous because
$T$ is. 
Altogether, 
$(x_n,u_n)_\nnn = (T(x_n+u_n),(\Id-T)(x_n+u_n))_\nnn$ converges weakly
to $(x,u) = (T(x+u),(\Id-T)(x+u))$. 
Thus, $(x,u)\in\gr A$ and therefore $\gr A$ is sequentially weakly closed. 
Conversely, let us assume that $\gr A$ is sequentially weakly closed.
Let $(x_n)_\nnn$ be a sequence in $X$ that is weakly convergent to $x$.
Our goal is to show that $Tx_n\weakly Tx$. 
Since $T$ is nonexpansive, the sequence $(Tx_n)_\nnn$ is bounded.
After passing to a subsequence and relabeling if necessary, 
we can and do assume that $(Tx_n)_\nnn$ converges weakly to some point
$y\in X$. 
Now $(Tx_n,x_n-Tx_n)_\nnn$ lies in $\gr A$, and this sequence
converges weakly to $(y,x-y)$. Since $\gr A$ is sequentially weakly closed,
it follows that $(y,x-y)\in\gr A$. 
Therefore, $x-y\in Ay$ 
$\Leftrightarrow$ $x\in (\Id+A)y$
$\Leftrightarrow$ $y=Tx$, which implies the result. 
\end{proof}

\begin{example}
\label{ex:skew}
Concerning items \ref{t:1viii} and \ref{t:1Banach}
in Theorem~\ref{t:1}, it is clear that
if $A$ is strongly monotone, then $T$ is a Banach contraction.
(This result is now part of the folklore; for perhaps the first
occurrence see \cite[p.~879f]{Rock76} and also 
\cite[Section~12.H]{Rock98} for further results in this direction.)
The converse, however, is false:
indeed, consider the case $X=\RR^2$ and set
\begin{equation}
\label{e:skew}
A = \begin{pmatrix}
0 & -1\\
1 & 0
\end{pmatrix}. 
\end{equation}
Then $(\forall z\in X)$ $\scal{z}{Az}=0$ so
$A$ cannot be strongly monotone.
On the other hand,
\begin{equation}
T = J_A = (\Id+A)^{-1} = 
\frac{1}{2}\begin{pmatrix}
1 & 1\\
-1 & 1
\end{pmatrix}
\end{equation}
is linear and $\|Tz\|^2 = \frac{1}{2}\|z\|^2$, which implies that
$T$ is a Banach contraction with constant $1/\sqrt{2}$. 
\end{example}

\begin{corollary}
\label{c:110113}
Let $A\colon X\to X$ be continuous, linear, and 
maximally monotone. 
Then the following hold.
\begin{enumerate}
\item 
\label{c:110113i}
If $J_A$ is a Banach contraction, then $A$ is (disjointly) injective.
\item
\label{c:110113ii}
If $\ran A$ is closed and $A$ is (disjointly) injective, 
then $J_A$ is a Banach contraction.
\end{enumerate}
\end{corollary}
\begin{proof}
The result is trivial if $X=\{0\}$ so we assume that $X\neq\{0\}$.

\ref{c:110113i}:
Assume that $J_A$ is a Banach contraction, say with constant
$\beta\in\left[0,1\right[$.
If $\beta = 0$, then $J_A \equiv 0$ $\Leftrightarrow$
$A=N_{\{0\}}$, which contradicts the single-valuedness of $A$. 
Thus, $0<\beta<1$. 
By Theorem~\ref{t:1}\ref{t:1Banach}, 
\begin{equation}
\label{e:110113:a}
(\forall x\in X)(\forall y\in X)
\quad
\frac{1-\beta^2}{\beta^2}\|x-y\|^2
\leq 2\scal{x-y}{Ax-Ay}+\|Ax-Ay\|^2.
\end{equation}
If $x\neq y$, then the left side of \eqref{e:110113:a}
is strictly positive, which implies that $Ax\neq Ay$.
Thus, $A$ is (disjointly) injective. 

\ref{c:110113ii}:
Let us assume that $\ran A$ is closed and that $A$ is 
(disjointly) injective.
Then $\ker A = \{0\}$ and hence, 
by e.g.\ \cite[Theorem~8.18]{Deutsch}, there exists
$\rho\in\RPP$ such that 
$(\forall z\in X)$
$\|Az\|\geq \rho\|z\|$.
Thus, 
$(\forall z\in X)$
$\|Az\|^2- \rho^2\|z\|\geq 0 $.
Set $\beta = 1/\sqrt{1+\rho^2}$. Then
$\rho^2 = (1-\beta^2)/\beta^2$ and hence
\begin{equation}
(\forall x\in X)(\forall y\in X)
\quad
\frac{1-\beta^2}{\beta^2}\|x-y\|^2
\leq \|Ax-Ay\|^2
\leq 2\scal{x-y}{Ax-Ay}+\|Ax-Ay\|^2.
\end{equation}
Again by Theorem~\ref{t:1}\ref{t:1Banach}, $J_A$ is a Banach
contraction with constant $\beta\in\zeroun$. 
\end{proof}

In item~\ref{c:110113ii} of Corollary~\ref{c:110113}, 
the assumption that $\ran A$ is closed is critical:

\begin{example}
Suppose that $X=\ell_2$, the space of
square-summable sequences, i.e., $x=(x_n)\in X$ if and only if
$\sum_{n=1}^\infty |x_n|^2 <\pinf$, and 
set 
\begin{equation}
A\colon X\to X\colon (x_n)\mapsto \big(\tfrac{1}{n}x_n\big).
\end{equation}
Then $A$ is continuous, linear, maximally monotone, and $\ran A$ is a 
dense, proper subspace of $X$ that is not closed. 
The resolvent $T=J_A$ is
\begin{equation}
T\colon X\to X\colon (x_n)\mapsto \big(\tfrac{n}{n+1}x_n\big).
\end{equation}
Now denote the $n^\mathrm{th}$ unit vector in $X$ by $\mathbf{e}_n$.
Then $\|T\mathbf{e}_n-T0\| = \tfrac{n}{n+1}\|\mathbf{e}_n-0\|$.
Since $\tfrac{n}{n+1}\to 1$, it follows that $T$ is not a Banach
contraction.
\end{example}

When $A$ is a subdifferential operator, then it is 
impossible to get the behavior witnessed in Example~\ref{ex:skew}:

\begin{proposition}
\label{p:india}
Let $f\in\Gamma_0$ and let $\varepsilon\in\RPP$.
Then $(1+\varepsilon)\prox_f$ is firmly nonexpansive
if and only if 
$\prox_f$ is a Banach contraction with constant $(1+\varepsilon)^{-1}$.
\end{proposition}
\begin{proof}
Set $\beta = (1+\varepsilon)^{-1}$.
It is clear that if  $(1+\varepsilon)\prox_f$ is firmly nonexpansive,
then  $(1+\varepsilon)\prox_f$ is nonexpansive and hence
$\prox_f$ is a Banach contraction with constant $\beta$. 
Conversely, assume that $\prox_f$ is a Banach contraction with constant
$\beta$. Since $\prox_f$
is the Fr\'echet gradient mapping of the continuous convex function
$f^*\Box\tfrac{1}{2}\|\cdot\|^2\colon X\to\RR$,
the Baillon-Haddad theorem (see \cite{BH} and also \cite{BC2010} for recent
results) guarantees that $\beta^{-1}\prox_f$ is firmly
nonexpansive. 
\end{proof}

\begin{remark}
If $n=2$, then
\eqref{e:t1cyc} reduces to
\begin{equation}
\scal{x_1-Tx_1}{Tx_1-Tx_2} + \scal{x_2-Tx_2}{Tx_2-Tx_1}\geq 0,
\end{equation}
i.e., to firm nonexpansiveness of $T$ (see
Fact~\ref{f:firm}\ref{f:firmsymm}). 
\end{remark}

\section{Duality and Self-Duality}

\label{s:3}

In the introduction, we pointed out that there is a natural duality
for firmly nonexpansive mappings and maximally monotone operators;
namely, $T\mapsto \Id-T$ and $A\mapsto A^{-1}$, respectively.
Thus, every property considered in Section~\ref{s:2} has a dual
property. We have considered all dual properties and
we shall explicitly single those out that we found to have 
simple and pleasing descriptions. 
Among these properties, those that are ``self-dual'', i.e.,
the property is identical to its dual property, stand
out even more. Let us begin my making these ideas more precise.

\begin{definition}[dual and self-dual properties]
Let $(p)$ and $(p^*)$ be properties for firmly 
nonexpansive mappings defined on $X$. 
If, for every firmly nonexpansive mapping $T\colon X\to X$,
\begin{equation}
\text{$T$ satisfies $(p)$ if and only if $\Id-T$ satisfies $(p^*)$,}
\end{equation}
then $(p^*)$ is \textbf{dual} to $(p)$, and hence $(p)$ is dual to $(p^*)$.
If $(p)=(p^*)$, we say that $(p)$ is \textbf{self-dual}.
Analogously, let $(q)$ and $(q^*)$ be properties of maximally monotone
operators defined on $X$.
If
\begin{equation}
\text{$A$ satisfies $(q)$ if and only if $A^{-1}$ satisfies $(q^*)$}
\end{equation}
for every maximally monotone operator $A\colon X\To X$, 
then $(q^*)$ is \textbf{dual} to $(q)$, and hence $(q)$ is dual to $(q^*)$.
If $(q)=(q^*)$, we say that $(q)$ is \textbf{self-dual}.
\end{definition}

We start with a simple example to illustrate these ideas.
It appears that here there is not natural dual
property to surjectivity of $T$ (except for the obvious
``$\Id-(\Id-T)$ is surjective'', but this is not illuminating).

\begin{theorem} \
\label{t:3}
Let $T\colon X\to X$ be firmly nonexpansive,
let $A\colon X\To X$ be maximally monotone,
and suppose that $T=J_A$ or equivalently that $A=T^{-1}-\Id$.
Then the following are equivalent:
\begin{enumerate}
\item
\label{t:3a}
$T$ is surjective.
\item
\label{t:3b}
$A$ has full domain.
\item
\label{t:3d}
$A^{-1}$ is surjective. 
\end{enumerate}
Thus for maximally monotone operators, 
surjectivity and full domain are properties that are dual to each
other. These properties are 
\textbf{not self-dual}:
e.g., $A=0$ has full domain while $A^{-1} = \partial\iota_{\{0\}}$ 
does not.
\end{theorem}
\begin{proof}
\ref{t:3a}$\Leftrightarrow$\ref{t:3b}: Theorem~\ref{t:1}\ref{t:1ii}.
\ref{t:3b}$\Leftrightarrow$\ref{t:3d}: Obvious.
\end{proof}

\begin{theorem} \ 
\label{t:2}
Let $T\colon X\to X$ be firmly nonexpansive,
let $A\colon X\To X$ be maximally monotone,
and suppose that $T=J_A$ or equivalently that $A=T^{-1}-\Id$.
Then the following are equivalent:
\begin{enumerate}
\item 
\label{t:2a}
$T$ is strictly nonexpansive.
\item 
\label{t:2b}
$A$ is disjointly injective.
\item 
\label{t:2c}
$\Id-T$ is injective.
\item 
\label{t:2d}
$A^{-1}$ is at most single-valued.
\end{enumerate}
Thus for firmly nonexpansive mappings,
strict nonexpansiveness and injectivity are dual to each other; 
and correspondingly for maximally monotone operators disjoint 
injectivity and at most single-valuedness are dual to each other.
These properties are \textbf{not self-dual}:
e.g., $T\equiv 0$ is strictly nonexpansive,
but $\Id-T = \Id$ is not;
correspondingly, $A = \partial\iota_{\{0\}}$ is 
disjointly injective but $A^{-1}=0$ is not.
\end{theorem}
\begin{proof}
We know that \ref{t:2a}$\Leftrightarrow$\ref{t:2b} by 
Theorem~\ref{t:1}\ref{t:1vi}.
We also know that \ref{t:2c}$\Leftrightarrow$\ref{t:2d} by 
Theorem~\ref{t:1}\ref{t:1iii} (applied to $A^{-1}$ and $\Id-T$).
It thus suffices to show that
\ref{t:2b}$\Leftrightarrow$\ref{t:2d}. 
Assume first that $A$ is disjointly injective and that
$\{x,y\}\subseteq A^{-1}u$. Then $u\in Ax\cap Ay$.
Since $A$ is disjointly injective, we have $x=y$.
Thus, $A^{-1}$ is at most single-valued.
Conversely, assume that $A^{-1}$ is at most single-valued and that
$u\in Ax\cap Ay$. Then $\{x,y\}\subseteq A^{-1}u$ and so $x=y$.
It follows that $A$ is disjointly injective.
\end{proof}

\begin{theorem}
\label{t:5}
Let $T\colon X\to X$ be firmly nonexpansive,
let $A\colon X\To X$ be maximally monotone,
and suppose that $T=J_A$ or equivalently that $A=T^{-1}-\Id$.
Then the following are equivalent:
\begin{enumerate}
\item
\label{t:5a}
$T$ satisfies \eqref{e:imagestrictfirm}.
\item
\label{t:5b}
$A$ is strictly monotone. 
\item
\label{t:5c}
$\Id-T$ satisfies
\begin{multline}
\label{e:coimsfirm}
(\forall x\in X)(\forall y\in X)\quad
\big(\Id-(\Id-T)\big)x\neq \big(\Id-(\Id-T)\big)y\\
\;\Rightarrow\;
\|(\Id-T)x-(\Id-T)y\|^2 < \scal{x-y}{(\Id-T)x-(\Id-T)y}.
\end{multline}
\item
\label{t:5d}
$A^{-1}$ satisfies
\begin{equation}
\label{e:costrictmonoA}
(\forall(x,u)\in\gr A^{-1})(\forall (y,v)\in\gr A^{-1})\quad
u\neq v \;\Rightarrow\;
\scal{x-y}{u-v}>0.
\end{equation}
\end{enumerate}
Thus for firmly nonexpansive mappings, properties 
\eqref{e:imagestrictfirm} and \eqref{e:coimsfirm}
are dual to each other;
and correspondingly for maximally monotone operators
strict monotonicity and \eqref{e:costrictmonoA} are dual to each
other. 
These properties are \textbf{not self-dual}:
e.g., $T = 0$ trivially satisfies \eqref{e:imagestrictfirm},
but $\Id-0 = \Id$ does not.
\end{theorem}
\begin{proof}
\ref{t:5a}$\Leftrightarrow$\ref{t:5b}: Theorem~\ref{t:1}\ref{t:1iv}. 
\ref{t:5a}$\Leftrightarrow$\ref{t:5c}: 
Indeed, \eqref{e:coimsfirm} and 
\eqref{e:imagestrictfirm} are equivalent as is easily seen by
expansion and rearranging.
\ref{t:5b}$\Leftrightarrow$\ref{t:5d}: 
Clear. 
\end{proof}

\begin{theorem}[self-duality of strict firm nonexpansiveness] 
\label{t:sundaya}
Let $T\colon X\to X$ be firmly nonexpansive,
let $A\colon X\To X$ be maximally monotone,
and suppose that $T=J_A$ or equivalently that $A=T^{-1}-\Id$.
Then the following are equivalent:
\begin{enumerate}
\item
$T$ is strictly firmly nonexpansive.
\item
$A$ is at most single-valued and strictly monotone.
\item
$\Id-T$ is strictly firmly nonexpansive.
\item
$A^{-1}$ is at most single-valued and strictly monotone.
\end{enumerate}
Consequently,
strict firm nonexpansive is a self-dual property for firmly
nonexpansive mappings; correspondingly, 
being both strictly monotone and  at most single-valued 
is self-dual for maximally monotone operators.
\end{theorem}
\begin{proof}
Note that 
$T$ is strictly firmly nonexpansive if and only if
\begin{equation}
(\forall x\in X)(\forall y\in X)\quad
x\neq y\;\Rightarrow\;
0<\scal{Tx-Ty}{(\Id-T)x-(\Id-T)y},
\end{equation}
which is obviously self-dual.
In view of Theorem~\ref{t:1}\ref{t:1v}, the corresponding property for $A$
is being both at most single-valued and strictly monotone.
\end{proof}

Theorem~\ref{t:sundaya} illustrates the technique of obtaining
self-dual properties by fusing any property and its dual.
Here is another example of this type.

\begin{theorem}[self-duality of strict nonexpansiveness and
injectivity]
\label{t:sundayb}
Let $T\colon X\to X$ be firmly nonexpansive,
let $A\colon X\To X$ be maximally monotone,
and suppose that $T=J_A$ or equivalently that $A=T^{-1}-\Id$.
Then the following are equivalent:
\begin{enumerate}
\item
$T$ is strictly nonexpansive and injective.
\item
$A$ is at most single-valued and disjointly injective.
\item
$\Id-T$ is strictly nonexpansive and injective.
\item
$A^{-1}$ is at most single-valued and disjointly injective.
\end{enumerate}
Consequently,
being both strictly nonexpansive and injective 
is a self-dual property for firmly
nonexpansive mappings; correspondingly, 
being both disjointly injective and  at most single-valued 
is self-dual for maximally monotone operators.
\end{theorem}
\begin{proof}
Clear from Theorem~\ref{t:2}\ref{t:1vii}. 
\end{proof}

\begin{remark} 
\label{r:sunday}
Some comments regarding Theorem~\ref{t:sundaya}
and Theorem~\ref{t:sundayb} are in order.
\begin{enumerate}
\item 
\label{r:sundayi}
Arguing directly (or by using the characterization with 
monotone operators via Theorem~\ref{t:1}), 
it is easy to verify the implication
\begin{equation}
\label{e:r:sunday}
\text{
$T$ is strictly firmly nonexpansive
$\Rightarrow$
$T$ is injective and strictly nonexpansive. }
\end{equation}
\item 
\label{r:sundayii}
The converse of implication \eqref{e:r:sunday} is false in
general: 
e.g., suppose that $X=\RR^2$, and 
let $A$ denote the counter-clockwise rotation by $\pi/2$, which we
utilized already in \eqref{e:skew}.
Clearly, $A$ is a linear single-valued maximally monotone operator that is
(disjointly) injective, but $A$ is not strictly monotone.
Accordingly, $T=J_A$ is linear, injective and strictly nonexpansive,
but not strictly firmly nonexpansive.
\item 
\label{r:sundayiii}
In striking contrast, we shall see 
shortly (in Corollary~\ref{c:toskana} 
below) that
when $X$ is finite-dimensional and $T=J_A$ is a {proximal mapping} 
(i.e., $A$ is a subdifferential operator), then 
the converse implication of \eqref{e:r:sunday} is true.
We shall write $\prox_f$ for $J_{\partial f}$ when
$f\in\Gamma_0$. 
\end{enumerate}
\end{remark}

\begin{lemma}
\label{l:toskana}
Suppose that $X$ is finite-dimensional and let $f\in\Gamma_0$.
Then the following are equivalent\footnote{See
\cite[Chapter~26]{Rock70} and \cite[Section~12.C]{Rock98}
for more on functions that are
essentially smooth, essentially strictly convex, or of Legendre type.}:
\begin{enumerate}
\item
\label{l:toskana:i}
$\partial f$ is disjointly injective.
\item
\label{l:toskana:ii}
$(\partial f)^{-1} = \partial f^*$ is at most single-valued.
\item
\label{l:toskana:iii}
$f^*$ is essentially smooth.
\item
\label{l:toskana:iv}
$f$ is essentially strictly convex.
\item
\label{l:toskana:v}
$\partial f$ is strictly monotone.
\item
\label{l:toskana:vi}
$\prox_f$ is strictly nonexpansive.
\item
\label{l:toskana:vii}
$(\forall x\in X)(\forall y\in X)$
$\prox_fx\neq \prox_fy$\\
$\Rightarrow$
$\|\prox_fx-\prox_fy\|^2 < \scal{x-y}{\prox_fx-\prox_fy}$.
\end{enumerate}
\end{lemma}
\begin{proof}
``\ref{l:toskana:i}$\Leftrightarrow$\ref{l:toskana:ii}'':
Theorem~\ref{t:2}.
``\ref{l:toskana:ii}$\Leftrightarrow$\ref{l:toskana:iii}'':
\cite[Theorem~26.1]{Rock70}. 
``\ref{l:toskana:iii}$\Leftrightarrow$\ref{l:toskana:iv}'':
\cite[Theorem~26.3]{Rock70}. 
``\ref{l:toskana:iv}$\Leftrightarrow$\ref{l:toskana:v}'':
\cite[Theorem~12.17]{Rock98}. 
``\ref{l:toskana:i}$\Leftrightarrow$\ref{l:toskana:vi}'':
Theorem~\ref{t:1}\ref{t:1vi}.
``\ref{l:toskana:v}$\Leftrightarrow$\ref{l:toskana:vii}'':
Theorem~\ref{t:1}\ref{t:1iv}.
\end{proof}

Lemma~\ref{l:toskana} admits a dual counterpart that 
contains 
various characterizations of essential smoothness.
The following consequence of these characterizations is also
related to Remark~\ref{r:sunday}.

\begin{corollary}[Legendre self-duality]
\label{c:toskana}
Suppose that $X$ is finite-dimensional and let
$f\in\Gamma_0$. 
Then the following are equivalent: 
\begin{enumerate}
\item
\label{c:toskana:i}
$\partial f$ is disjointly injective and at most single-valued.
\item
\label{c:toskana:ii}
$\partial f$ is strictly monotone and at most single-valued.
\item
\label{c:toskana:iii}
$f$ is Legendre.
\item
\label{c:toskana:iv}
$\prox_f$ is strictly firmly nonexpansive.
\item
\label{c:toskana:v}
$\prox_f$ is strictly nonexpansive and injective.
\item
\label{c:toskana:vi}
$\partial f^*$ is disjointly injective and at most single-valued.
\item
\label{c:toskana:vii}
$\partial f^*$ is strictly monotone and at most single-valued.
\item
\label{c:toskana:viii}
$f^*$ is Legendre.
\item
\label{c:toskana:ix}
$\prox_{f^*}$ is strictly firmly nonexpansive.
\item
\label{c:toskana:x}
$\prox_{f^*}$ is strictly nonexpansive and injective.
\end{enumerate}
\end{corollary}

\begin{remark}
When $X$ is infinite-dimensional, the results become more technical 
and additional assumptions are required 
due to subtleties of Legendre functions that do not occur 
in finite-dimensional settings; 
see \cite{BBC} and \cite{BorVanBook}.
\end{remark}

\begin{remark}[strong monotonicity and cocoercivity] 
Concerning Theorem~\ref{t:1}\ref{t:1viii}\&\ref{t:1coco},
we note in passing that strong monotonicity is not self-dual:
indeed, $A = \partial\iota_{\{0\}}$ is strongly monotone, but
$A^{-1}\equiv 0$ is not. This example also shows that uniform
monotonicity is not self-dual. The property dual to strong monotonicity is
cocoercivity (which is also known the more accurate name 
\textbf{inverse strong monotonicity}).
\end{remark}

\begin{theorem}[self-duality of paramonotonicity]
\label{t:para}
Let $A\colon X\To X$ be maximally monotone,
let $T\colon X\to X$ be firmly nonexpansive, and
suppose that $T=J_A$ or equivalently that $A=T^{-1}-\Id$.
Then $A$ is paramonotone if and only if $A^{-1}$;
consequently, $T$ satisfies \eqref{e:fpara} if and only if $\Id-T$
satisfies \eqref{e:fpara} (with $T$ replaced by $\Id-T$).
Consequently,
being paramonotone is a
is a self-dual property for maximally monotone operators;
correspondingly, satisfying \eqref{e:fpara} is a self-dual property
for firmly nonexpansive mappings.
\end{theorem}
\begin{proof}
Self-duality is immediate from the definition of paramonotonicity,
and the corresponding result for firmly nonexpansive mappings follows
from Theorem~\ref{t:1}\ref{t:1para}.
\end{proof}

\begin{theorem}[self-duality of cyclical firm nonexpansiveness and
cyclical monotonicity] \ \\
\label{t:subdiff}
Let $T\colon X\to X$ be firmly nonexpansive,
let $A\colon X\To X$ be maximally monotone, 
let $f\in\Gamma_0$, 
and suppose that
$T=J_A$ or equivalently that $A=T^{-1}-\Id$.
Then the following are equivalent:
\begin{enumerate}
\item $T$ is cyclically firmly nonexpansive.
\item $A$ is cyclically monotone.
\item $A=\partial f$.
\item $\Id-T$ is cyclically firmly nonexpansive.
\item $A^{-1}$ is cyclically monotone.
\item $A^{-1}=\partial f^*$.
\end{enumerate}
Consequently,
cyclic firm nonexpansiveness is a self-dual property for
firmly nonexpansive mappings; 
correspondingly, cyclic monotonicity is
is a self-dual property for maximally monotone operators.
\end{theorem}
\begin{proof}
The fact that cyclically maximal monotone operators are
subdifferential operators is due to Rockafellar and well known,
as is the identity $(\partial f)^{-1}=\partial f^*$.
The result thus follows from Theorem~\ref{t:1}\ref{t:1cyc}. 
\end{proof}

\begin{theorem}[self-duality of rectangularity]
\label{t:3*}
Let $T\colon X\to X$ be firmly nonexpansive,
let $A\colon X\To X$ be maximally monotone, 
and suppose that
$T=J_A$ or equivalently that $A=T^{-1}-\Id$.
Then the following are equivalent:
\begin{enumerate}
\item $T$ satisfies \eqref{e:t13*}.
\item $A$ is rectangular.
\item $\Id-T$ satisfies \eqref{e:t13*}.
\item $A^{-1}$ is rectangular.
\end{enumerate}
Consequently, rectangularity is a self-dual property for
maximally monotone operators;
correspondingly, \eqref{e:t13*} is a self-dual property for
firmly nonexpansive mappings. 
\end{theorem}
\begin{proof}
It is obvious from the definition that either property is self-dual;
the equivalences thus follow from Theorem~\ref{t:1}\ref{t:13*}.
\end{proof}

\begin{theorem}[self-duality of linearity]
\label{t:linear}
Let $T\colon X\to X$ be firmly nonexpansive,
let $A\colon X\To X$ be maximally monotone, 
and suppose that
$T=J_A$ or equivalently that $A=T^{-1}-\Id$.
Then the following are equivalent:
\begin{enumerate}
\item $T$ is linear.
\item $A$ is a linear relation. 
\item $\Id-T$ is linear.
\item $A^{-1}$ is a linear relation. 
\end{enumerate}
Consequently,
linearity is a self-dual property for
firmly nonexpansive mappings; 
correspondingly, being a linear relation is
a self-dual property for maximally monotone operators.
\end{theorem}
\begin{proof}
It is clear that $T$ is linear if and only if $\Id-T$ is;
thus, the result follows from 
Theorem~\ref{t:1}\ref{t:1linear}.
\end{proof}

\begin{theorem}[self-duality of affineness]
\label{t:affine}
Let $T\colon X\to X$ be firmly nonexpansive,
let $A\colon X\To X$ be maximally monotone, 
and suppose that
$T=J_A$ or equivalently that $A=T^{-1}-\Id$.
Then the following are equivalent:
\begin{enumerate}
\item $T$ is affine.
\item $A$ is an affine relation. 
\item $\Id-T$ is affine.
\item $A^{-1}$ is an affine relation. 
\end{enumerate}
Consequently,
affineness is a self-dual property for
firmly nonexpansive mappings; 
correspondingly, being an affine relation is
is a self-dual property for maximally monotone operators.
\end{theorem}
\begin{proof}
It is clear that $T$ is affine if and only if $\Id-T$ is;
therefore, the result follows from 
Theorem~\ref{t:1}\ref{t:1affine}.
\end{proof}

\begin{remark}[projection]
Concerning Theorem~\ref{t:1}\ref{t:1zara}, we note in passing
that being a projection is not a self-dual: indeed, 
suppose that $X\neq\{0\}$ and 
let $T$ be the projection onto the closed unit ball. 
Then $\Id-T$ is not a
projection since $\Fix(\Id-T)=\{0\}\subsetneqq X = \ran(\Id-T)$. 
\end{remark}

\begin{theorem}[self-duality of sequential weak continuity]
\label{t:weak}
Let $T\colon X\to X$ be firmly nonexpansive,
let $A\colon X\To X$ be maximally monotone, 
and suppose that
$T=J_A$ or equivalently that $A=T^{-1}-\Id$.
Then the following are equivalent:
\begin{enumerate}
\item $T$ is sequentially weakly continuous.
\item $\gr A$ is sequentially weakly closed.
\item $\Id-T$ is sequentially weakly continuous.
\item $A^{-1}$ is sequentially weakly closed. 
\end{enumerate}
Consequently,
sequential weak continuity is a self-dual property for
firmly nonexpansive mappings; 
correspondingly, having a sequentially weakly closed graph
is a self-dual property for maximally monotone operators.
\end{theorem}
\begin{proof}
Since $\Id$ is weakly continuous, 
it is clear that $T$ is sequentially weakly continuous if and only if
$\Id-T$ is; thus, the result follows from
Theorem~\ref{t:1}\ref{t:1weak}. 
\end{proof}

We end this section by listing all self-dual properties encountered so
far. 

\begin{remark}
Let $T\colon X\to X$ be firmly nonexpansive,
let $A\colon X\To X$ be maximally monotone, 
and suppose that
$T=J_A$ or equivalently that $A=T^{-1}-\Id$.
Then the following properties are self-dual:
\begin{enumerate}
\item
$T$ is strictly firmly nonexpansive (Theorem~\ref{t:sundaya}). 
\item 
$T$ is strictly nonexpansive and injective (Theorem~\ref{t:sundayb}).
\item 
$A$ is paramonotone (Theorem~\ref{t:para}).
\item 
$A$ is a subdifferential operator (Theorem~\ref{t:subdiff}).
\item
$A$ is rectangular (Theorem~\ref{t:3*}). 
\item
$T$ is linear (Theorem~\ref{t:linear}).
\item
$T$ is affine (Theorem~\ref{t:affine}).
\item
$T$ is sequentially weakly continuous (Theorem~\ref{t:weak}). 
\end{enumerate}
\end{remark}

\section{Nonexpansive Mappings}

\label{s:4}

In the previous two sections, 
we have extensively utilized the the correspondence between
firmly nonexpansive mappings and maximally monotone operators.
However, Fact~\ref{f:firm} provides another correspondence, namely
with nonexpansive mappings.
Indeed,
\begin{equation}
\text{
$T$ is firmly nonexpansive if and only if 
$N=2T-\Id$ is nonexpansive.
}
\end{equation}
Note that $N$ is also referred to as a \textbf{reflected resolvent}. 
The corresponding dual of $N$ within the set of nonexpansive mappings is,
very elegantly, simply  
\begin{equation}
-N.
\end{equation}
Thus, all results have counterparts formulated for nonexpansive
mappings. These counterparts are most easily derived from
the firmly nonexpansive formulation, by simply replacing
$T$ by $\tfrac{1}{2}\Id+\tfrac{1}{2}N$.
This makes the proofs fairly straightforward and we accordingly omit
them. 
In the following, we list those cases where the 
corresponding nonexpansive formulation turned out to be
simple and elegant.

\begin{theorem}[strict firm nonexpansiveness]
\label{t:reisera}
Let $T\colon X\to X$ be firmly nonexpansive,
let $N \colon X\to X$ be nonexpansive,
and suppose that $N=2T-\Id$. 
Then $T$ is strictly firmly nonexpansive if and only if 
$N$ is strictly nonexpansive.
\end{theorem}

\begin{remark} \ 
\begin{enumerate}
\item
Of course, from Theorem~\ref{t:sundaya}, we know that 
strict firm nonexpansiveness is a self-dual property. This can also be
seen within the realm of nonexpansive mappings since $N$ is strictly
nonexpansive if and only if $-N$ is.
\item 
Furthermore, combining Theorem~\ref{t:sundaya} with
Theorem~\ref{t:reisera} yields the following:
a maximally monotone operator $A$ is at most single-valued and
strictly monotone if and only if its reflected resolvent $2J_A-\Id$ is
strictly nonexpansive.
This characterization was observed by Rockafellar and Wets;
see \cite[Proposition~12.11]{Rock98}. 
\item In passing, we note that when $X$ is finite-dimensional,
the iterates of a strictly nonexpansive
mapping converge to the unique fixed point (assuming it exists).
For this and more, see, e.g., \cite{EKN}. 
\end{enumerate}
\end{remark}

The next result pertains to Theorem~\ref{t:1}\ref{t:1viii}. 

\begin{theorem}[strong monotonicity]
\label{t:reiserb}
Let $A\colon X\To X$ be maximally monotone, 
let $N \colon X\to X$ be nonexpansive,
suppose that $N=2J_A-\Id$ and that $\varepsilon\in\RPP$.
Then $A$ is strongly monotone with constant $\varepsilon$
if and only if $\varepsilon\Id + (1+\varepsilon)N$ is nonexpansive.
\end{theorem}

\begin{theorem}[reflected resolvent as Banach contraction]
\label{t:Banachrr}
Let $A\colon X\To X$ be maximally monotone, 
let $T\colon X\to X$ be firmly nonexpansive,
and let $N \colon X\to X$ be nonexpansive.
Suppose that $T=J_A$, that $N=2T-\Id$, and
that $\beta\in\left[0,1\right]$. 
Then the following are equivalent:
\begin{enumerate}
\item
\label{t:Banachrri}
$(\forall (x,u)\in\gr A)(\forall (y,v)\in\gr A)$\\
$(1-\beta^2)(\|x-y\|^2 + \|u-v\|^2)
\leq 2(1+\beta^2)\scal{x-y}{u-v}$
\item
\label{t:Banachrrii}
$(\forall x\in X)(\forall y\in X)$
$(1-\beta^2)\|x-y\|^2 \leq 4\scal{Tx-Ty}{(\Id-T)x-(\Id-T)y}$
\item
\label{t:Banachrriii}
$(\forall x\in X)(\forall y\in X)$
$\|Nx-Ny\|\leq\beta\|x-y\|$.
\end{enumerate}
\end{theorem}
\begin{proof}
In view of the Minty parametrization (see \eqref{e:Mintypar}),
item~\ref{t:Banachrri} is equivalent to 
\begin{multline}
\label{e:Banachrr}
(\forall x\in X)(\forall y\in X)
\quad
(1-\beta^2)
\big(\|Tx-Ty\|^2 + \|(x-Tx)-(y-Ty)\|^2\big)\\
\leq 2\big(1+\beta^2\big)\scal{Tx-Ty}{(x-Tx)-(y-Ty)}.
\end{multline}
Simple algebraic manipulations show that 
\eqref{e:Banachrr} is equivalent to \ref{t:Banachrrii},
which in turn is equivalent to \ref{t:Banachrriii}. 
\end{proof}

It is clear that the properties \ref{t:Banachrri}--\ref{t:Banachrriii}
in Theorem~\ref{t:Banachrr} are self-dual (for fixed $\beta$).
The following result is a simple consequence.

\begin{corollary}[self-duality of reflected resolvents that are Banach
contractions] 
\label{c:Banachrr}
Let $A\colon X\To X$ be maximally monotone, 
let $T\colon X\to X$ be firmly nonexpansive,
let $N \colon X\to X$ be nonexpansive,
and suppose that $T=J_A$ and $N=2T-\Id$.
Then the following are equivalent:
\begin{enumerate}
\item
\label{c:Banachrri}
$\displaystyle \inf
\mmenge{ \frac{
\scal{x-y}{u-v}}{\|x-y\|^2+\|u-v\|^2}}{\{(x,u),(y,v)\}\subseteq\gr A
\text{ and } (x,u)\neq(y,v)} > 0$. 
\item
$\displaystyle \inf
\mmenge{ \frac{
\scal{Tx-Ty}{(\Id-T)x-(\Id-T)y}}{\|x-y\|^2}}{\{x,y\}\subseteq X
\text{ and } x\neq y} > 0$. 
\item $N$ is a Banach contraction. 
\end{enumerate}
Furthermore, these properties are self-dual for their respective
classes of operators. 
\end{corollary}

\begin{remark} 
Precisely when  $A\colon x\mapsto x-z$ for some fixed vector
$z\in X$, we compute $T \colon x\mapsto (x+z)/2$ and therefore
we reach the extreme case of Corollary~\ref{c:Banachrr} where
$N\colon x\mapsto z$ is a Banach contraction with constant $0$.
\end{remark}

\begin{corollary}
\label{c:monday}
Let $A\colon X\To X$ be maximally monotone, 
let $T\colon X\to X$ be firmly nonexpansive,
let $N \colon X\to X$ be nonexpansive,
and suppose that $T=J_A$ and $N=2T-\Id$.
Then the following are equivalent:
\begin{enumerate}
\item 
\label{c:mondayi}
Both $A$ and $A^{-1}$ are strongly monotone.
\item 
\label{c:mondayii}
There exists $\gamma\in\left]1,\pinf\right[$ such that both 
$\gamma T$ and $\gamma(\Id-T)$ are firmly nonexpansive.
\item 
\label{c:mondayiii}
$N$ is a Banach contraction.
\end{enumerate}
\end{corollary}
\begin{proof}
Let us assume that $A$ and $A^{-1}$ are both strongly monotone;
equivalently, there exists $\varepsilon\in\RPP$ such that
$A-\varepsilon\Id$ and $A^{-1}-\varepsilon\Id$ are monotone. 
Let $\{(x,u),(y,v)\}\subseteq\gr A$.
Then $\{(u,x),(v,y)\}\subseteq\gr A^{-1}$ and
\begin{equation}
\scal{x-y}{u-v} \geq \varepsilon\|x-y\|^2
\;\text{and}\;
\scal{u-v}{x-y} \geq \varepsilon \|u-v\|^2.
\end{equation}
Adding these inequalities yields $2\scal{x-y}{u-v} \geq
\varepsilon(\|x-y\|^2+\|u-v\|^2)$. 
Thus, item~\ref{c:Banachrri} of Corollary~\ref{c:Banachrr} holds.
Conversely, if item~\ref{c:Banachrri} of Corollary~\ref{c:Banachrr} holds,
then both $A$ and $A^{-1}$ are strongly monotone. 
Therefore, by Corollary~\ref{c:Banachrr}, \ref{c:mondayi} and
\ref{c:mondayiii} are equivalent. 
Finally, in view of Theorem~\ref{t:1}\ref{t:1viii}, we see that
\ref{c:mondayi} and \ref{c:mondayii} are also equivalent. 
\end{proof}

Additional characterizations are available for subdifferential operators:

\begin{proposition}
\label{p:tuesday}
Let $f\in\Gamma_0$. Then the following are equivalent:
\begin{enumerate}
\item
\label{p:tuesdayi}
$f$ and $f^*$ are strongly convex.
\item
\label{p:tuesdayii}
$f$ and $f^*$ are everywhere differentiable, and both $\nabla f$ and
$\nabla f^*$ are Lipschitz continuous.
\item
\label{p:tuesdayiii}
$\prox_f$ and $\Id-\prox_f$ are Banach contractions.
\item
\label{p:tuesdayiv}
$2\prox_f-\Id$ is a Banach contraction.
\end{enumerate}
\end{proposition}
\begin{proof}
It is well known that for functions, strong convexity is equivalent to
strong monotonicity of the subdifferential operators; see, e.g., 
\cite[Example~22.3]{BC2011}. 
In view of Proposition~\ref{p:india} and Corollary~\ref{c:monday},
we obtain the equivalence of items \ref{p:tuesdayi},
\ref{p:tuesdayiii}, and \ref{p:tuesdayiv}. 
Finally, the equivalence of \ref{p:tuesdayi} and \ref{p:tuesdayii}
follows from \cite[Theorem~2.1]{BC2010}. 
\end{proof}

We now turn to linear relations.

\begin{proposition}
\label{p:linrel}
Let $A\colon X\to X$ be a maximally monotone linear relation.
Then the following are equivalent:
\begin{enumerate}
\item
\label{p:linreli}
Both $A$ and $A^{-1}$ are strongly monotone.
\item 
\label{p:linrelii}
$A$ is a continuous surjective linear operator on $X$ and 
$\displaystyle 
\inf_{z\in X\smallsetminus \{0\}} \frac{\scal{z}{Az}}{\|z\|^2 + \|Az\|^2} > 0$.
\item
\label{p:linreliii}
$2J_A-\Id$ is a Banach contraction.
\end{enumerate}
If $X$ is finite-dimensional, then {\rm \ref{p:linreli}--\ref{p:linreliii}}
are also equivalent to 
\begin{enumerate}
\setcounter{enumi}{3}
\item 
\label{p:linreliv}
$A\colon X\to X$ satisfies $(\forall z\in X\smallsetminus\{0\})$
$\scal{z}{Az}>0$. 
\end{enumerate}
\end{proposition}
\begin{proof}
``\ref{p:linreli}$\Leftrightarrow$\ref{p:linreliii}'':  
Clear from Corollary~\ref{c:monday}.

``\ref{p:linreli}$\Rightarrow$\ref{p:linrelii}'':
By, e.g.\ \cite[Example~22.9(iii)]{BC2011}, $A$ and $A^{-1}$ are single-valued
surjective operators with full domain. Since $A$ and $A^{-1}$ are linear,
\cite[Corollary~21.19]{BC2011} implies that $A$ and $A^{-1}$ are
continuous. Thus, \ref{p:linrelii} holds. 

``\ref{p:linreli}$\Leftarrow$\ref{p:linrelii}'':
\ref{p:linrelii} implies that item~\ref{c:Banachrri} of
Corollary~\ref{c:Banachrr} holds. Thus, \ref{p:linreli} follows from
Corollary~\ref{c:Banachrr} and Corollary~\ref{c:monday}.

``\ref{p:linrelii}$\Rightarrow$\ref{p:linreliv}'': Clear.

``\ref{p:linrelii}$\Leftarrow$\ref{p:linreliv}'': 
Since $A$ is injective and $X$ is finite-dimensional, $A$ is bijective and
continuous. To see that the infimum in item~\ref{p:linrelii} is strictly
positive, note that we may take the infimum over the unit sphere, which is
a compact subset of $X$. 
\end{proof}

We shall conclude this paper with some comments regarding
prototypical applications of the above results to splitting methods 
(see also \cite{BC2011} for further information and various variants). 
Here is a technical lemma, which is part of the folklore,
and whose simple proof we omit.

\begin{lemma}
\label{l:110115:1}
Let $T_1,\ldots,T_n$ be finitely many 
nonexpansive mappings from $X$ to $X$,
and let $\lambda_1,\ldots,\lambda_n$ be in $\left]0,1\right]$ such that
$\lambda_1+\cdots + \lambda_n=1$.
Then the following hold:
\begin{enumerate}
\item
The composition
$T_1T_2\cdots T_n$ is nonexpansive.
\item 
The convex combination 
$\lambda_1T_1+\cdots+\lambda_nT_n$ is nonexpansive.
\item
\label{l:110115:1i}
If some $T_i$ is strictly nonexpansive, then 
$T_1T_2\cdots T_n$ is strictly nonexpansive.
\item
\label{l:110115:1ii}
If some $T_i$ is strictly nonexpansive, then 
$\lambda_1T_1+\cdots+\lambda_nT_n$ is strictly nonexpansive.
\item
If some $T_i$ is a Banach contraction, then 
$T_1T_2\cdots T_n$ is a Banach contraction. 
\item
If some $T_i$ is a Banach contraction, then 
$\lambda_1T_1+\cdots+\lambda_nT_n$ is a Banach contraction.
\end{enumerate}
\end{lemma}

\begin{corollary}[backward-backward iteration]
Let $A_1$ and $A_2$ be two 
maximally monotone operators from $X$ to $X$,
and assume that 
one of these is disjointly injective. 
Then the
(backward-backward) composition $T_1T_2$ is strictly nonexpansive.
\end{corollary}
\begin{proof}
Combine 
Theorem~\ref{t:1}\ref{t:1vi} and Lemma~\ref{l:110115:1}. 
\end{proof}

\begin{corollary}[Douglas-Rachford iteration]
\label{c:DR}
Let $A_1$ and $A_2$ be two maximally monotone operators from $X$ to
$X$, 
and assume that one of these is both at most single-valued
and strictly monotone (as is, e.g., the subdifferential operator of
a Legendre function; see Corollary~\ref{c:toskana}).
Denote the resolvents of $A_1$ and $A_2$ by $T_1$ and $T_2$,
respectively. 
Then the operator governing the Douglas-Rachford
iteration, i.e., 
\begin{equation}
\label{e:DR}
T := \tfrac{1}{2}(2T_1-\Id)(2T_2-\Id) + \tfrac{1}{2}\Id, 
\end{equation}
is not just firmly nonexpansive but also strictly
nonexpansive; consequently, $\Fix T$ is either empty or a singleton.
\end{corollary}
\begin{proof}
In view of Theorem~\ref{t:sundaya} and Theorem~\ref{t:reisera}, 
we see that
$2T_1-\Id$ and $2T_2-\Id$ are both nonexpansive, and one of these two
is strictly nonexpansive.
By Lemma~\ref{l:110115:1}\ref{l:110115:1i},
$(2T_1-\Id)(2T_2-\Id)$ is strictly nonexpansive.
Hence, by Lemma~\ref{l:110115:1}\ref{l:110115:1ii},
$T$ is strictly nonexpansive. 
\end{proof}

\begin{remark} 
Consider Corollary~\ref{c:DR}, and assume that $A_i$, where
$i\in\{1,2\}$, 
satisfies condition~\ref{c:Banachrri} in Corollary~\ref{c:Banachrr}.
Then $2T_i-\Id$ is a Banach contraction by Corollary~\ref{c:Banachrr}.
Furthermore, Lemma~\ref{l:110115:1} now shows that 
the Douglas-Rachford operator $T$ defined in \eqref{e:DR}
is a Banach contraction. Thus, $\Fix T$ is a singleton and the unique
fixed point may be found as the strong limit of any
sequence of Banach-Picard iterates for $T$. 
See also \cite[Section~25.2]{BC2011} for various variants and
strengthenings. 
\end{remark}

\section*{Acknowledgments}
Heinz Bauschke was partially supported by the Natural Sciences and
Engineering Research Council of Canada and by the Canada Research Chair
Program.
Xianfu Wang was partially
supported by the Natural Sciences and Engineering Research Council
of Canada.


\small 
\begin{thebibliography}{999}

\bibitem{BH}
J.-B.\ Baillon and G.\ Haddad,
Quelques propri\'et\'es des op\'erateurs angle-born\'es et
$n$-cycliquement monotones,
\emph{Israel Journal of Mathematics} 26 (1977), 137--150. 

\bibitem{BBBRW}
S.\ Bartz, H.H.\ Bauschke, J.M.\ Borwein,
S.\ Reich, and X.\ Wang,
Fitzpatrick functions, cyclic monotonicity and Rockafellar's
antiderivative,
\emph{Nonlinear Analysis} 66 (2007), 1198--1223.

\bibitem{BBC}
H.H.\ Bauschke, J.M.\ Borwein, and P.L.\ Combettes,
Essential smoothness, essential strict convexity,
and Legendre functions in Banach spaces,
\emph{Communications in Contemporary Mathematics} 3 (2001), 615--647. 

\bibitem{BC2010}
H.H.\ Bauschke and P.L.\ Combettes,
The Baillon-Haddad theorem revisited,
\emph{Journal of Convex Analysis} 17 (2010), 781--787.

\bibitem{BC2011}
H.H.\ Bauschke and P.L.\ Combettes,
\emph{Convex Analysis and Monotone Operator Theory in Hilbert Spaces},
Springer-Verlag, 2011. 


\bibitem{BorVanBook}
J.M.\ Borwein and J.D.\ Vanderwerff,
\emph{Convex Functions}, 
Cambridge University Press, 2010.

\bibitem{Brezis}
H. Br\'ezis,
\emph{Operateurs Maximaux Monotones et
Semi-Groupes de Contractions dans les Espaces de Hilbert},
North-Holland/Elsevier, 1973. 

\bibitem{Browder67}
F.E.\ Browder,
Convergence theorems for sequences of nonlinear operators
in Banach spaces,
\emph{Mathematische Zeitschrift} 100 (1967), 201--225. 

\bibitem{BurIus}
R.S.\ Burachik and A.N.\ Iusem,
\emph{Set-Valued Mappings and Enlargements
of Monotone Operators},
Springer-Verlag, 2008. 

\bibitem{CIZ}
Y. Censor, A.N.\ Iusem, and S.A.\ Zenios,
An interior point method with Bregman functions
for the variational inequality problem with 
paramonotone operators,
\emph{Mathematical Programming Series~A}~81 (1998), 373--400.

\bibitem{Cross}
R.\ Cross,
\emph{Multivalued Linear Operators},
Marcel Dekker, 1998.

\bibitem{Deutsch}
F.\ Deutsch,
\emph{Best Approximation in Inner Product Spaces},
Springer-Verlag, 2001.

\bibitem{EckBer}
J.\ Eckstein and D.P.\ Bertsekas,
\emph{On the Douglas-Rachford splitting method
and the proximal point algorithm for maximal monotone
operators},
\emph{Mathematical Programming Series A} 55 (1992), 293--318.

\bibitem{EKN}
L.\ Elsner, I.\ Koltracht, and M.\ Neumann,
Convergence of sequential and asynchronous
nonlinear paracontractions,
\emph{Numerische Mathematik}~62 (1992), 305--319.

\bibitem{GK}
K.\ Goebel and W.A.\ Kirk,
\emph{Topics in Metric Fixed Point Theory},
Cambridge University Press, 1990. 

\bibitem{GR}
K.\ Goebel and S.\ Reich,
\emph{Uniform Convexity, Hyperbolic Geometry, and Nonexpansive Mappings},
Marcel Dekker, 1984.

\bibitem{Minty}
G.J.\ Minty,
Monotone (nonlinear) operators in Hilbert spaces,
\emph{Duke Mathematical Journal} 29 (1962), 341--346.


\bibitem{Moreau}
J.-J.\ Moreau, 
Proximit\'e et dualit\'e dans un espace hilbertien,
\emph{Bulletin de la Soci\'et\'e Math\'ematique de France} 93 (1965),
273--299.

\bibitem{Rock70}
R.T.\ Rockafellar,
\emph{Convex Analysis},
Princeton University Press, Princeton, 1970.

\bibitem{Rock76}
R.T.\ Rockafellar,
Monotone operators and the proximal point algorithm,
\emph{SIAM Journal on Control and Optimization} 14
(1976), 877--898. 

\bibitem{Rock98}
R.T.\ Rockafellar and R. J-B\ Wets,
\emph{Variational Analysis},
Springer-Verlag, 
1998.

\bibitem{Simons1}
S.\ Simons,
\emph{Minimax and Monotonicity},
Springer-Verlag,
1998.

\bibitem{Simons2}
S.\ Simons,
\emph{From Hahn-Banach to Monotonicity},
Springer-Verlag,
2008.

\bibitem{Zalinescu}{C.\ Z\u{a}linescu},
\emph{Convex Analysis in General Vector Spaces},
World Scientific Publishing, 2002.

\bibitem{Zara}
E.H.\ Zarantonello,
Projections on convex sets in Hilbert space and spectral theory I.
Projections on convex sets, in
\emph{Contributions to Nonlinear Functional Analysis},
E.H.\ Zarantonello (editor), pp.~237--341, Academic Press, 1971. 


\bibitem{Zeidler2a}
E.\ Zeidler,
\emph{Nonlinear Functional Analysis and Its Applications II/A:
Linear Monotone Operators},
Springer-Verlag, 1990.

\bibitem{Zeidler2b}
E.\ Zeidler,
\emph{Nonlinear Functional Analysis and Its Applications II/B:
Nonlinear Monotone Operators},
Springer-Verlag, 1990.

\bibitem{Zeidler1}
E.\ Zeidler,
\emph{Nonlinear Functional Analysis and Its Applications I:
Fixed Point Theorems},
Springer-Verlag, 1993.

\end{thebibliography}
\end{document}